\documentclass{article}%
\usepackage{graphicx}
\usepackage[intlimits]{amsmath}
\usepackage{amsfonts}
\usepackage{amssymb}%
\newtheorem{corollary}{Corollary}

\newtheorem{theorem}{Theorem}

\newenvironment{proof}{\noindent{\bf Proof\,}}{\hfill$\Box$}

\begin{document}

\title{The sum of degrees in cliques}
\author{B\'{e}la Bollob\'{a}s\thanks{Department of Mathematical Sciences, University
of Memphis, Memphis TN 38152, USA} \thanks{Trinity College, Cambridge CB2 1TQ,
UK} \thanks{Research supported in part by DARPA grant F33615-01-C-1900.} \ and
Vladimir Nikiforov$^{\ast}$}
\maketitle

\begin{abstract}
For every graph $G,$ let%
\[
\Delta_{r}\left(  G\right)  =\max\left\{  \sum_{u\in R}d\left(  u\right)
:R\text{ is an }r\text{-clique of }G\right\}
\]
and let $\Delta_{r}\left(  n,m\right)  $ be the minimum of $\Delta_{r}\left(
G\right)  $ taken over all graphs of order $n$ and size $m$. Write
$t_{r}\left(  n\right)  $ for the size of the $r$-chromatic Tur\'{a}n graph of
order $n$.

Improving earlier results of Edwards and Faudree, we show that for every
$r\geq2,$ if $m\geq t_{r}\left(  n\right)  ,$ then%
\begin{equation}
\Delta_{r}\left(  n,m\right)  \geq\frac{2rm}{n}, \label{mineq}%
\end{equation}
as conjectured by Bollob\'{a}s and Erd\H{o}s.

It is known that inequality (\ref{mineq}) fails for $m<t_{r}\left(  n\right)
.$ However, we show that for every $\varepsilon>0,$ there is $\delta>0$ such
that if $m>t_{r}\left(  n\right)  -\delta n^{2}$ then
\[
\Delta_{r}\left(  n,m\right)  \geq\left(  1-\varepsilon\right)  \frac{2rm}%
{n}.
\]
Finally, we generalize (\ref{mineq}) to graphs with edge weights.

\end{abstract}

\section{Introduction}

Our notation and terminology are standard (see, e.g. \cite{Bol98}): thus
$G\left(  n,m\right)  $ stands for a graph of $n$ vertices and $m$ edges. For
a graph $G$ and a vertex $u\in V\left(  G\right)  ,$ we write $\Gamma\left(
u\right)  $ for the set of vertices adjacent to $u$ and set $d_{G}\left(
u\right)  =\left\vert \Gamma\left(  u\right)  \right\vert ;$ we write
$d\left(  u\right)  $ instead of $d_{G}\left(  u\right)  $ if the graph $G$ is
understood. However, somewhat unusually, for $U\subset V\left(  G\right)  ,$
we set $\widehat{\Gamma}\left(  U\right)  =\left\vert \cap_{v\in U}%
\Gamma\left(  v\right)  \right\vert $ and $\widehat{d}\left(  U\right)
=\left\vert \widehat{\Gamma}\left(  U\right)  \right\vert $.

We write $T_{r}\left(  n\right)  $ for the $r$-chromatic Tur\'{a}n graph on
$n$ vertices and $t_{r}\left(  n\right)  $ for the number of its edges.

For every $r\geq2$ and every graph $G,$ let $\Delta_{r}\left(  G\right)  $ be
the maximum of the sum of degrees of the vertices of an $r$-clique, as in the
Abstract. If $G$ has no $r$-cliques, we set $\Delta_{r}\left(  G\right)  =0.$
Furthermore, let
\[
\Delta_{r}\left(  n,m\right)  =\min_{G=G\left(  n,m\right)  }\Delta_{r}\left(
G\right)  .
\]

Since $T_{r}\left(  n\right)  $ is a $K_{r+1}$-free graph, it follows that
$\Delta_{r}\left(  n,m\right)  =0$ for $m\leq t_{r-1}\left(  n\right)  .$ In
1975 Bollob\'{a}s and Erd\H{o}s \cite{BoEr75} conjectured that for every
$r\geq2,$ if $m\geq t_{r}\left(  n\right)  ,$ then%
\begin{equation}
\Delta_{r}\left(  n,m\right)  \geq\frac{2rm}{n}. \label{conjBE}%
\end{equation}

Edwards \cite{Edw77}, \cite{Edw78} proved (\ref{conjBE}) under the weaker
condition $m>\left(  r-1\right)  n^{2}/2r;$ he also proved that the conjecture
holds for $2\leq r\leq8$ and $n\geq r^{2}$. Later Faudree \cite{Fdr92} proved
the conjecture for any $r\geq2$ and $n>r^{2}\left(  r-1\right)  /4$.

For $t_{r-1}\left(  n\right)  <m<t_{r}\left(  n\right)  $ the value of
$\Delta_{r}\left(  n,m\right)  $ is essentially unknown even for $r=3$ (see
\cite{ErLa83}, \cite{Fan88} and \cite{Fdr92} for partial results.) An explicit
construction due to Erd\H{o}s (see \cite{Fdr92}) shows that, for every
$\varepsilon>0,$ there exists $\delta>0$ such that if $t_{r-1}\left(
n\right)  <m<t_{r}\left(  n\right)  -\delta n^{2}$ then%
\[
\Delta_{r}\left(  n,m\right)  \leq\left(  1-\varepsilon\right)  \frac{2rm}%
{n}.
\]

In this note we prove a stronger form of (\ref{conjBE}) for every $r$ and $n.$
Furthermore, we prove that $\Delta_{r}\left(  n,m\right)  $ is
\textquotedblleft stable\textquotedblright\ as $m$ approaches $t_{r}\left(
n\right)  .$ More precisely, for every $\varepsilon>0,$ there is $\delta>0$
such that if $m>t_{r}\left(  n\right)  -\delta n^{2}$ then
\[
\Delta_{r}\left(  n,m\right)  \geq\left(  1-\varepsilon\right)  \frac{2rm}{n}%
\]
for $n$ sufficiently large.

\subsection{Preliminary observations}

Let $M_{1},...,M_{k}$ are subsets of a (finite) set $V$ with complements
$\overline{M}_{1},...,\overline{M}_{k}.$ Then
\[
\sum_{i=1}^{k}\left\vert \overline{M_{i}}\right\vert \geq\left\vert \cup
_{i=1}^{k}\overline{M_{i}}\right\vert
\]
and so,%
\begin{equation}
\left\vert \cap_{i=1}^{k}M_{i}\right\vert \geq\sum_{i=1}^{k}\left\vert
M_{i}\right\vert -\left(  k-1\right)  \left\vert V\right\vert . \label{basin}%
\end{equation}

The size $t_{r}\left(  n\right)  $ of the Tur\'{a}n graph $T_{r}\left(
n\right)  $ is given by%
\[
t_{r}\left(  n\right)  =\frac{r-1}{2r}\left(  n^{2}-s^{2}\right)  +\binom
{s}{2}.
\]
where $s$ is the reminder of $n$ modulo $r:$ Hence,
\begin{equation}
\frac{r-1}{2r}n^{2}\geq t_{r}\left(  n\right)  \geq\frac{r-1}{2r}n^{2}%
-\frac{r}{8}. \label{turest}%
\end{equation}

\section{A greedy algorithm}

In what follows we shall identify a clique with its vertex set.

Faudree \cite{Fdr92} introduced the following algorithm $\mathfrak{P}$ which
constructs a clique $\left\{  v_{1},...,v_{k}\right\}  $ in a graph $G$:

Step 1: $v_{1}$ is a vertex of maximum degree in $G$;

Step 2: having selected $v_{1},...,v_{i-1},$ if $\widehat{\Gamma}\left(
v_{1}...v_{i-1}\right)  =\varnothing$ then $\mathfrak{P}$ stops, otherwise
$\mathfrak{P}$ selects a vertex of maximum degree $v_{i}\in\widehat{\Gamma
}\left(  v_{1}...v_{i-1}\right)  $ and step 2 is repeated again.

Faudree's main reason to introduce this algorithm was to mainly to prove
Conjecture (\ref{conjBE}) for $n$ sufficiently large, so he did not study
$\mathfrak{P}$ in great detail. In this section we shall establish some
properties of $\mathfrak{P}$ for their own sake. Later, in Section \ref{deg},
we shall apply these results to prove an extension of (\ref{conjBE}) for every
$n.$

Note that $\mathfrak{P}$ need not construct a unique sequence. Sequences that
can be constructed by $\mathfrak{P}$ are called $\mathfrak{P}$%
\emph{-sequences; }the definition of $\mathfrak{P}$ implies that
$\widehat{\Gamma}\left(  v_{1}...v_{k}\right)  =\varnothing$ for every
$\mathfrak{P}$-sequence\emph{ }$v_{1},...,v_{k}.$

\begin{theorem}
\label{Turext} Let $r\geq2,$ $n\geq r$ and $m\geq t_{r}\left(  n\right)  $.
Then every graph $G=G\left(  n,m\right)  $ is such that:

\emph{(i)} every $\mathfrak{P}$-sequence has at least $r$ terms;

\emph{(ii)} for every $\mathfrak{P}$-sequence $v_{1},...,v_{r},$%
\begin{equation}
\sum_{i=1}^{r}d\left(  v_{i}\right)  \geq\left(  r-1\right)  n; \label{con0}%
\end{equation}

\emph{(iii)} if equality holds in (\ref{con0}) for some $\mathfrak{P}%
$-sequence $v_{1},...,v_{r}$ then $m=t_{r}\left(  n\right)  $.
\end{theorem}

\begin{proof}
Without loss of generality we may assume that $\mathfrak{P}$ constructs
exactly the vertices $1,...,k$ and hence $d\left(  1\right)  \geq...\geq
d\left(  k\right)  $.

\emph{Proof} $\emph{of}$ $\emph{(i)}$ $\emph{and}$ $\emph{(ii)}$ To prove
\emph{(i)} we have to show that $k\geq r.$ For every $i=1,...,k,$ let
$M_{i}=\Gamma\left(  i\right)  ;$ clearly,
\[
\sum_{i=1}^{k}d\left(  i\right)  \leq\left(  q-1\right)  n,
\]
since, otherwise, (\ref{basin}) implies that $\widehat{\Gamma}\left(
v_{1}...v_{k}\right)  \neq\varnothing,$ and so $1,...,k$ is not a
$\mathfrak{P}$-sequence, contradicting the choice of $k$. Suppose $k<r$, and
let $q$ be the smallest integer such that the inequality
\begin{equation}
\sum_{i=1}^{h}d\left(  i\right)  >\left(  h-1\right)  n \label{con1}%
\end{equation}
holds for $h=1,...,q-1,$ while%
\begin{equation}
\sum_{i=1}^{q}d\left(  i\right)  \leq\left(  q-1\right)  n. \label{con2}%
\end{equation}
Clearly, $1<q\leq k$.

Partition $V=\cup_{i=1}^{q}V_{i},$ so that
\begin{align*}
V_{1}  &  =V\backslash\Gamma\left(  1\right)  ,\\
\text{ }V_{i}  &  =\widehat{\Gamma}\left(  \left[  i-1\right]  \right)
\backslash\widehat{\Gamma}\left(  \left[  i\right]  \right)  \text{ \ \ for
\ }i=2,...,q-1,\\
\text{ }V_{q}  &  =\widehat{\Gamma}\left(  \left[  q-1\right]  \right)  .
\end{align*}

We have
\begin{align}
2m  &  =\sum_{j\in V}d\left(  j\right)  =\sum_{h=1}^{q}\sum_{j\in V_{h}%
}d\left(  j\right)  \leq\sum_{i=1}^{q}d\left(  i\right)  \left\vert
V_{i}\right\vert \nonumber\\
&  =d\left(  1\right)  \left(  n-d\left(  1\right)  \right)  +\sum_{i=2}%
^{q-1}d\left(  i\right)  \left(  \widehat{d}\left(  \left[  i-1\right]
\right)  -\widehat{d}\left(  \left[  i\right]  \right)  \right)  +d\left(
q\right)  \widehat{d}\left(  \left[  q-1\right]  \right) \nonumber\\
&  =d\left(  1\right)  n+\sum_{i=1}^{q-1}\widehat{d}\left(  \left[  i\right]
\right)  \left(  d\left(  i+1\right)  -d\left(  i\right)  \right)  .
\label{con3}%
\end{align}

For every $h\in\left[  q-2\right]  ,$ applying (\ref{basin}) with
$M_{i}=\Gamma\left(  i\right)  $, $i\in\left[  h\right]  ,$ we see that,%
\[
\widehat{d}\left(  \left[  h\right]  \right)  =\left\vert \widehat{\Gamma
}\left(  \left[  h\right]  \right)  \right\vert \geq\sum_{i=1}^{h}d\left(
i\right)  -\left(  h-1\right)  n=n-\sum_{i=1}^{h}\left(  n-d\left(  i\right)
\right)  >0,
\]
and hence, by $d\left(  h+1\right)  \leq d\left(  h\right)  ,$ it follows
that
\begin{equation}
\widehat{d}\left(  \left[  h\right]  \right)  \left(  d\left(  h+1\right)
-d\left(  h\right)  \right)  \leq\left(  n-\sum_{i=1}^{h}\left(  n-d\left(
i\right)  \right)  \right)  \left(  d\left(  h+1\right)  -d\left(  h\right)
\right)  . \label{con3.1}%
\end{equation}

Since, from (\ref{con2}), we have%
\begin{equation}
d\left(  q\right)  \leq\left(  q-1\right)  n-\sum_{i=1}^{q-1}d\left(
i\right)  =\sum_{i=1}^{q-1}\left(  n-d\left(  i\right)  \right)  ,
\label{con3.2}%
\end{equation}
in view of (\ref{basin}), (\ref{con1}), and (\ref{con3}), we deduce
\begin{align*}
\widehat{d}\left(  \left[  q-1\right]  \right)  \left(  d\left(  q\right)
-d\left(  q-1\right)  \right)   &  \leq\left(  n-\sum_{i=1}^{q-1}\left(
n-d\left(  i\right)  \right)  \right)  \left(  d\left(  q\right)  -d\left(
q-1\right)  \right) \\
&  \leq\left(  n-\sum_{i=1}^{q-1}\left(  n-d\left(  i\right)  \right)
\right)  \left(  \sum_{i=1}^{q-1}\left(  n-d\left(  i\right)  \right)
-d\left(  q-1\right)  \right)  .
\end{align*}

Recalling (\ref{con3.1}) and (\ref{con3.2}), this inequality implies that
\begin{align*}
2m  &  \leq nd\left(  1\right)  +\sum_{h=1}^{q-2}\left(  n-\sum_{i=1}%
^{h}\left(  n-d\left(  i\right)  \right)  \right)  \left(  d\left(
h+1\right)  -d\left(  h\right)  \right) \\
&  +\left(  n-\sum_{i=1}^{q-1}\left(  n-d\left(  i\right)  \right)  \right)
\left(  \sum_{i=1}^{q-1}\left(  n-d\left(  i\right)  \right)  -d\left(
q-1\right)  \right)  .
\end{align*}
Dividing by $2$ and rearranging the right-hand side, we obtain
\begin{equation}
m\leq\left(  n-\sum_{i=1}^{q-1}\left(  n-d\left(  i\right)  \right)  \right)
\left(  \sum_{i=1}^{q-1}\left(  n-d\left(  i\right)  \right)  \right)
+\sum_{1\leq i<j\leq q-1}\left(  n-d\left(  i\right)  \right)  \left(
n-d\left(  j\right)  \right)  . \label{con4}%
\end{equation}

Furthermore, for every $i\in\left[  q-1\right]  ,$ set $k_{i}=n-d\left(
i\right)  ;$ let $k_{q}=n-\left(  k_{1}+...+k_{q-1}\right)  .$ Clearly,
$k_{i}>0$ for every $i\in\left[  q\right]  $; also, $k_{1}+...+k_{q}=n$. Note
that the right-hand side of (\ref{con4}) is exactly%
\[
\sum_{1\leq i<j\leq q}k_{i}k_{j},
\]
and this is precisely $e\left(  K\left(  k_{1},...,k_{q}\right)  \right)  .$
Given $n$ and $k_{1}+...+k_{q}=n,$ the value $e\left(  K\left(  k_{1}%
,...,k_{q}\right)  \right)  $ attains its maximum if and only if all $k_{i}$
differ by at most $1,$ that is to say, when $K\left(  k_{1},...,k_{q}\right)
$ is exactly the Tur\'{a}n graph $T_{q}\left(  n\right)  .$ Hence, $m\geq
t_{r}\left(  n\right)  $ and (\ref{con4}) imply%
\begin{equation}
t_{r}\left(  n\right)  \leq m\leq e\left(  K\left(  k_{1},...,k_{q}\right)
\right)  \leq t_{q}\left(  n\right)  . \label{con4.1}%
\end{equation}
Since $q<r\leq n$ implies $t_{q}\left(  n\right)  <t_{r}\left(  n\right)  ,$
contradicting (\ref{con4.1}), the proof of \emph{(i)} is complete.

To prove \emph{(ii)} suppose (\ref{con0}) fails, i.e.,%
\[
\sum_{i=1}^{r}d\left(  i\right)  <\left(  r-1\right)  n.
\]
Hence, (\ref{con3.2}) holds with a strict inequality and so, the proof of
(\ref{con4.1}) gives $t_{r}\left(  n\right)  <t_{r}\left(  n\right)  .$ This
contradiction completes the proof of \emph{(ii)}.

\emph{Proof of (iii) }Suppose that for some $\mathfrak{P}$-sequence
$v_{1},...,v_{r}$ equality holds in (\ref{con0}). We may and shall assume that
$v_{1},...,v_{r}=1,...,r$, i.e.,%
\[
\sum_{i=1}^{r}d\left(  i\right)  =\left(  r-1\right)  n.
\]
Following the arguments in the proof of \emph{(i)} and \emph{(ii)}, from
(\ref{con4.1}) we conclude that%
\[
t_{r}\left(  n\right)  \leq m\leq t_{r}\left(  n\right)  .
\]
and this completes the proof.
\end{proof}

\section{\label{deg}Degree sums in cliques}

In this section we turn to the problem of finding $\Delta_{r}\left(
n,m\right)  $ for $m\geq t_{r}\left(  n\right)  .$ We shall apply Theorem
\ref{Turext} to prove that every graph $G=G\left(  n,m\right)  $ with $m\geq
t_{r}\left(  n\right)  $ contains an $r$-clique $R$ with%
\begin{equation}
\sum_{i\in R}d\left(  i\right)  \geq\frac{2rm}{n}. \label{BEin}%
\end{equation}
As proved by Faudree \cite{Fdr92}, the required $r$-clique $R$ may be
constructed by the algorithm $\mathfrak{P}$. Note that the assertion is
trivial for regular graphs; as we shall show, if $G$ is not regular, we may
demand strict inequality in (\ref{BEin}).

\begin{theorem}
\label{Edwext} Let $r\geq2,$ $n\geq r,$ $m\geq t_{r}\left(  n\right)  $ and
let $G=G\left(  n,m\right)  $ be a graph which is not regular. Then there
exists a $\mathfrak{P}$-sequence $v_{1},...,v_{r}$ such that%
\[
\sum_{i=1}^{r}d\left(  v_{i}\right)  >\frac{2rm}{n}.
\]

\end{theorem}

\begin{proof}
Part \emph{(iii) }of Theorem \ref{Turext} implies that for some $\mathfrak{P}%
$-sequence of $r$ vertices, say $1,...,r,$ we have%
\[
\sum_{i=1}^{r}d\left(  i\right)  >\left(  r-1\right)  n.
\]
Since $d\left(  i\right)  <n,$ we immediately obtain%
\begin{equation}
\sum_{i=1}^{s}d\left(  i\right)  >\left(  s-1\right)  n \label{cond1}%
\end{equation}
for every $s\in\left[  r\right]  .$

The rest of the proof consists of two parts: In part \emph{(a)} we find an
upper bound for $m$ in terms of $\sum_{i=1}^{r}d\left(  i\right)  $ and
$\sum_{i=1}^{r}d^{2}\left(  i\right)  .$ Then, in part \emph{(b)}, we prove
that%
\[
\frac{1}{r}\sum_{i=1}^{r}d\left(  i\right)  \geq\frac{2m}{n},
\]
and show that if equality holds then $G$ is regular.

\emph{(a)} Partition the set $V$ into $r$ sets $V=V_{1}\cup...\cup V_{r},$
where,
\begin{align*}
V_{1}  &  =V\backslash\Gamma\left(  1\right)  ,\\
V_{i}  &  =\widehat{\Gamma}\left(  \left[  i-1\right]  \right)  \backslash
\widehat{\Gamma}\left(  \left[  i\right]  \right)  \text{ for }%
i=2,..,r-2,\text{ }\\
V_{r}  &  =\widehat{\Gamma}\left(  \left[  r-1\right]  \right)  .
\end{align*}
We have,
\begin{align}
2m  &  =\sum_{i\in V}d\left(  i\right)  =\sum_{h=1}^{r}\sum_{j\in V_{h}%
}d\left(  j\right)  \leq\sum_{i=1}^{r}d\left(  i\right)  \left\vert
V_{i}\right\vert \nonumber\\
&  =\sum_{i=1}^{r-1}\left(  d\left(  i\right)  -d\left(  r\right)  \right)
\left\vert V_{i}\right\vert +nd\left(  r\right)  \label{cond3}%
\end{align}
Clearly, for every $i\in\left[  r-1\right]  ,$ from (\ref{basin}), we have
\[
\left\vert \widehat{\Gamma}\left(  \left[  i+1\right]  \right)  \right\vert
\geq\left\vert \widehat{\Gamma}\left(  \left[  i\right]  \right)  \right\vert
+\left\vert \Gamma\left(  i+1\right)  \right\vert -n=\left\vert \widehat
{\Gamma}\left(  \left[  i\right]  \right)  \right\vert +d\left(  i+1\right)
-n
\]
and hence, $\left\vert V_{i}\right\vert \leq n-d\left(  i\right)  $ holds for
every $i\in\left[  r-1\right]  .$ Estimating $\left\vert V_{i}\right\vert $ in
(\ref{cond3}) we obtain%
\begin{align*}
2m  &  \leq\sum_{i=1}^{r-1}\left(  d\left(  i\right)  -d\left(  r\right)
\right)  \left(  n-d\left(  i\right)  \right)  +nd\left(  r\right) \\
&  =n\sum_{i=1}^{r}d\left(  i\right)  -\sum_{i=1}^{r}d^{2}\left(  i\right)
+d\left(  r\right)  \left(  \sum_{i=1}^{r}d\left(  i\right)  -n\left(
r-1\right)  \right)  .
\end{align*}

\emph{(b) }Let $S_{r}=\sum_{i=1}^{r}d\left(  i\right)  .$ From $d\left(
r\right)  \leq S_{r}/r$ and Cauchy's inequality we deduce%
\begin{align*}
2m  &  \leq nS_{r}-\sum_{i=1}^{r}d^{2}\left(  i\right)  +\frac{S_{r}}%
{r}\left(  S_{r}-\left(  r-1\right)  n\right) \\
&  \leq nS_{r}-\frac{1}{r}\left(  S_{r}\right)  ^{2}+\frac{S_{r}}{r}\left(
S_{r}-\left(  r-1\right)  n\right)  \leq\frac{nS_{r}}{r},
\end{align*}
and so,%
\begin{equation}
\sum_{i=1}^{r}d\left(  i\right)  \geq\frac{2rm}{n}. \label{BEin2}%
\end{equation}

To complete the proof suppose we have an equality in (\ref{BEin2}). This
implies that%
\[
\sum_{i=1}^{r}d^{2}\left(  i\right)  =\frac{1}{r}\left(  \sum_{i=1}%
^{r}d\left(  i\right)  \right)  ^{2}%
\]
and so, $d\left(  1\right)  =...=d\left(  r\right)  .$ Therefore, the maximum
degree $d\left(  1\right)  $ equals the average degree $2m/n,$ contradicting
the assumption that $G$ is not regular.
\end{proof}

Since for every $m\geq t_{r}\left(  n\right)  $ there is a graph $G=G\left(
n,m\right)  $ whose degrees differ by at most 1, we obtain the following
bounds on $\Delta_{r}\left(  n,m\right)  .$

\begin{corollary}
For every $m\geq t_{r}\left(  n\right)  $%
\[
\frac{2rm}{n}\leq\Delta_{r}\left(  n,m\right)  <\frac{2rm}{n}+r.
\]

\end{corollary}

\section{Stability of $\Delta_{r}\left(  n,m\right)  $ as $m$ approaches
$t_{r}\left(  n\right)  $}

It is known that inequality (\ref{conjBE}) is far from being true \ if $m\leq
t_{r}\left(  n\right)  -\varepsilon n$ for some $\varepsilon>0$ (e.g., see
\cite{Fdr92}). However, it turns out that, as $m$ approaches $t_{r}\left(
n\right)  ,$ the function $\Delta_{r}\left(  n,m\right)  $ approaches $2rm/n.$
More precisely, the following stability result holds.

\begin{theorem}
\label{Tstab}For every $\varepsilon>0$ there exist $n_{0}=n_{0}\left(
\varepsilon\right)  $ and $\delta=\delta\left(  \varepsilon\right)  >0$ such
that if $m>t_{r}\left(  n\right)  -\delta n^{2}$ then
\[
\Delta_{r}\left(  n,m\right)  >\left(  1-\varepsilon\right)  \frac{2rm}{n}%
\]
for all $n>n_{0}.$
\end{theorem}

\begin{proof}
Without loss of generality we may assume that%
\[
0<\varepsilon<\frac{2}{r\left(  r+1\right)  }.
\]
Set
\[
\delta=\delta\left(  \varepsilon\right)  =\frac{1}{32}\varepsilon^{2}.
\]

If $m\geq t_{r}\left(  n\right)  ,$ the assertion follows from Theorem
\ref{Edwext}, hence we may assume that
\[
\frac{2rm}{n}<\frac{2rt_{r}\left(  n\right)  }{n}\leq\left(  r-1\right)  n.
\]

Clearly, our theorem follows if we show that $m>t_{r}\left(  n\right)  -\delta
n^{2}$ implies
\begin{equation}
\Delta_{r}\left(  n,m\right)  >\left(  1-\varepsilon\right)  \left(
r-1\right)  n \label{ldeg}%
\end{equation}
for $n$ sufficiently large.

Suppose the graph $G=G\left(  n,m\right)  $ satisfies $m>t_{r}\left(
n\right)  -\delta n^{2}.$\ By (\ref{turest}), if $n$ is large enough,
\begin{equation}
m>t_{r}\left(  n\right)  -\delta n^{2}>\left(  \frac{r-1}{2r}-\delta\right)
n^{2}-\frac{r}{8}\geq\left(  \frac{r-1}{2r}-2\delta\right)  n^{2}.
\label{mbnd}%
\end{equation}

Let $M_{\varepsilon}\subset V$ be defined as
\[
M_{\varepsilon}=\left\{  u:d\left(  u\right)  \leq\left(  \frac{r-1}{r}%
-\frac{\varepsilon}{2}\right)  n\right\}  .
\]

The rest of the proof consists of two parts. In part \emph{(a) }we shall show
that $\left\vert M_{\varepsilon}\right\vert <\varepsilon n,$ and in part
\emph{(b)} we shall show that the subgraph induced by $V\backslash
M_{\varepsilon}$ contains an $r$-clique with large degree sum, proving
(\ref{ldeg}).

\emph{(a) }Our first goal is to show that $\left\vert M_{\varepsilon
}\right\vert <\varepsilon n.$ Indeed, assume the opposite and select an
arbitrary $M^{\prime}\subset M_{\varepsilon}$ satisfying
\begin{equation}
\left(  \frac{1}{2}-\frac{1}{2\sqrt{2}}\right)  \varepsilon n<\left\vert
M^{\prime}\right\vert <\left(  \frac{1}{2}+\frac{1}{2\sqrt{2}}\right)
\varepsilon n. \label{bnds}%
\end{equation}
Let $G^{\prime}$ be the subgraph of $G$ induced by $V\backslash M^{\prime}.$
Then%
\begin{align}
e\left(  G\right)   &  =e\left(  G^{\prime}\right)  +e\left(  M^{\prime
},V\backslash M^{\prime}\right)  +e\left(  M^{\prime}\right)  \leq e\left(
G^{\prime}\right)  +\sum_{u\in M^{\prime}}d\left(  u\right) \label{ebnd}\\
&  \leq e\left(  G^{\prime}\right)  +\left\vert M^{\prime}\right\vert \left(
\frac{r-1}{r}-\frac{\varepsilon}{2}\right)  n.\nonumber
\end{align}

Observe that second inequality of (\ref{bnds}) implies%
\[
n-\left\vert M^{\prime}\right\vert >\left(  1-\varepsilon\right)  n.
\]

Hence, if
\[
e\left(  G^{\prime}\right)  \geq\frac{r-1}{2r}\left(  n-\left\vert M^{\prime
}\right\vert \right)  ^{2}%
\]
then, applying Theorem \ref{Edwext} to the graph $G^{\prime}$, we see that
\[
\Delta_{r}\left(  G\right)  \geq\Delta_{r}\left(  G^{\prime}\right)  \geq
\frac{2re\left(  G^{\prime}\right)  }{n-\left\vert M^{\prime}\right\vert }%
\geq\left(  r-1\right)  \left(  n-\left\vert M^{\prime}\right\vert \right)
>\left(  r-1\right)  \left(  1-\varepsilon\right)  n,
\]
and (\ref{ldeg}) follows. Therefore, we may assume
\[
e\left(  G^{\prime}\right)  <\frac{r-1}{2r}\left(  n-\left\vert M^{\prime
}\right\vert \right)  ^{2}.
\]
Then, by (\ref{mbnd}) and (\ref{ebnd}),
\[
\frac{r-1}{2r}\left(  n-\left\vert M^{\prime}\right\vert \right)
^{2}>e\left(  G^{\prime}\right)  >-\left\vert M^{\prime}\right\vert \left(
\frac{r-1}{r}-\frac{\varepsilon}{2}\right)  n+\left(  \frac{r-1}{2r}%
-2\delta\right)  n^{2}.
\]

Setting $x=\left\vert M^{\prime}\right\vert /n,$ this shows that%
\[
\frac{r-1}{2r}\left(  1-x\right)  ^{2}+x\left(  \frac{r-1}{r}-\frac
{\varepsilon}{2}\right)  -\left(  \frac{r-1}{2r}-2\delta\right)  >0,
\]
which imply that%
\[
x^{2}-\varepsilon x+4\delta>0.
\]
Hence, either%
\[
\left\vert M^{\prime}\right\vert >\left(  \frac{\varepsilon-\sqrt
{\varepsilon^{2}-16\delta}}{2}\right)  n=\left(  \frac{1}{2}-\frac{1}%
{2\sqrt{2}}\right)  \varepsilon n
\]
or%
\[
\left\vert M^{\prime}\right\vert <\left(  \frac{\varepsilon+\sqrt
{\varepsilon^{2}-16\delta}}{2}\right)  =\left(  \frac{1}{2}+\frac{1}{2\sqrt
{2}}\right)  \varepsilon n,
\]
contradicting (\ref{bnds}). Therefore, $\left\vert M_{\varepsilon}\right\vert
<\varepsilon n,$ as claimed

\emph{(b)} Let $G_{0}$ be the subgraph of $G$ induced by $V\backslash
M_{\varepsilon}.$ By the definition of $M_{\varepsilon},$ if $u\in V\backslash
M_{\varepsilon},$ then%
\[
d_{G}\left(  u\right)  >\left(  \frac{r-1}{r}-\frac{\varepsilon}{2}\right)
n,
\]
and so%
\[
d_{G_{0}}\left(  u\right)  >\left(  \frac{r-1}{r}-\frac{\varepsilon}%
{2}\right)  n-\left\vert M_{\varepsilon}\right\vert >\frac{r-2}{r-1}\left(
n-\left\vert M_{\varepsilon}\right\vert \right)  .
\]
Hence, by Tur\'{a}n's theorem, $G_{0}$ contains an $r$-clique and, therefore,%
\[
\Delta_{r}\left(  G\right)  >r\left(  \frac{r-1}{r}-\frac{\varepsilon}%
{2}\right)  n\geq\left(  1-\varepsilon\right)  \left(  r-1\right)  n,
\]
proving (\ref{ldeg}) and completing the proof of our theorem.
\end{proof}

\textbf{Acknowledgement} The authors thank Prof. D. Todorov for pointing out a
fallacy in an earlier version of the proof of Theorem \ref{Edwext}.

\end{document}